\documentclass[a4paper,12pt]{amsart}
\usepackage{}
\usepackage{ifthen}
\usepackage{amsfonts}
\usepackage{amssymb}
\usepackage{amsmath,amssymb,graphicx}
\usepackage{color}
\numberwithin{equation}{section}
\nonstopmode
\newtheorem{thm}{Theorem}[section]

\newtheorem{lem}[thm]{Lemma}

\newtheorem{cor}[thm]{Corollary}

\newcommand{\be}{\begin{eqnarray}}
\newcommand{\ee}{\end{eqnarray}}
\newcommand{\ba}{\begin{array}}
\newcommand{\ea}{\end{array}}
\newcommand{\ben}{\begin{eqnarray*}}
\newcommand{\een}{\end{eqnarray*}}

\newcommand{\pa}{\partial}

\newcommand{\lam}{\lambda}

\renewcommand{\arg}{{\operatorname{arg}\,}}

\newcommand{\C}{{\mathbb C}}
\newcommand{\D}{{\mathbb D}}

\newcommand{\uhp}{{\mathbb H}}

\newcommand{\R}{{\mathbb R}}

\newcommand{\Z}{{\mathbb Z}}

\newcommand{\bD}{{\overline{\mathbb D}}}
\newcommand{\sphere}{{\widehat{\mathbb C}}}

\renewcommand{\Im}{\,{\operatorname{Im}\,}}

\renewcommand{\Re}{{\operatorname{Re}\,}}

\renewcommand{\mod}{{\operatorname{mod}\,}}
\newcommand{\inv}{^{-1}}

\newcommand{\dist}{{\operatorname{dist}}}
\newcommand{\arth}{{\,\operatorname{arth}\,}}
\newcommand{\arctanh}{{\operatorname{arth}\,}}

\newcommand{\PSL}{{\operatorname{PSL}}}
\newcommand{\SL}{{\operatorname{SL}}}

\newcommand{\aand}{{\quad\text{and}\quad}}

\newenvironment{pf}[1][]{%
 \vskip 3mm
 \noindent
 \ifthenelse{\equal{#1}{}}%
  {{\slshape Proof. }}%
  {{\slshape #1.} }%
 }%
{\qed\bigskip}

\newcounter{minutes}
\divide\time by 60
\newcounter{hours}
\multiply\time by 60
\addtocounter{minutes}{-\time}
\begin{document}
\bibliographystyle{amsplain}
\title{
Construction of nearly hyperbolic distance on punctured spheres
}

%\begin{center}
%{\tiny \texttt{FILE:~\jobname .tex,
%        printed: \number\year-\number\month-\number\day,
%        \thehours.\ifnum\theminutes<10{0}\fi\theminutes}
%}
%\end{center}

\author[T. Sugawa]{Toshiyuki Sugawa}
\address{Graduate School of Information Sciences,
Tohoku University, Aoba-ku, Sendai 980-8579, Japan}
\email{sugawa@math.is.tohoku.ac.jp}
\author[T. Zhang]{Tanran Zhang}
\address{Department of Mathematics, Soochow University, No.1 Shizi Street, Suzhou 215006, China}
\email{trzhang@suda.edu.cn}
\keywords{hyperbolic metric, punctured sphere}
\subjclass[2010]{Primary 30C35; Secondary 30C55}
\begin{abstract}
We define a distance function on the bordered punctured disk $0<|z|\le 1/e$
in the complex plane, which is comparable with the hyperbolic distance
of the punctured unit disk $0<|z|<1.$
As an application, we will construct a distance function on an $n$-times
punctured sphere which is comparable with the hyperbolic distance.
We also propose a comparable quantity which is not necessarily a distance function on the
punctured sphere but easier to compute.
\end{abstract}
\thanks{
The authors were supported in part by JSPS Grant-in-Aid for
Scientific Research (B) 22340025.
}
\maketitle

\section{Introduction}

A domain $\Omega$ in the Riemann sphere $\sphere=\C\cup\{\infty\}$ has
the upper half-plane $\uhp=\{\zeta\in\C: \Im\zeta>0\}$
as its holomorphic universal covering space
precisely if the complement $\sphere\setminus\Omega$ contains at least three points.
Such a domain is called {\it hyperbolic}.
Since the Poincar\'e metric $|d\zeta|/(2\Im \zeta)$ is invariant under the
pullback by analytic automorphisms of $\uhp$ (in other words,
under the action of $\PSL(2,\R)$), it descends to a metric on $\Omega,$
called the {\it hyperbolic metric} of $\Omega$ and denoted by
$\lambda_\Omega(w)|dw|.$
More explicitly, they are related by the formula
$\lambda_\Omega(p(\zeta))|p'(\zeta)|=1/(2\Im \zeta),$ where
$p:\uhp\to\Omega$ is a holomorphic universal covering projection
of $\uhp$ onto $\Omega.$
The quantity $\lambda_\Omega(w)$ is sometimes called the {\it hyperbolic
density} of $\Omega$ and it is independent of the particular choice of
$p$ and $\zeta\in p\inv(w).$
Note that $\lambda_\Omega$ has constant Gaussian curvature $-4$ on $\Omega.$
We denote by $h_\Omega(w_1,w_2)$ the distance function induced by
$\lambda_\Omega,$
called the {\it hyperbolic distance} on $\Omega$ and the distance is known to be complete.
That is, $h_\Omega(w_1,w_2)=\inf_\alpha \ell_\Omega(\alpha),$
where the infimum is taken over all rectifiable curves $\alpha$
joining $w_1$ and $w_2$ in $\Omega$ and
$$
\ell_\Omega(\alpha)=\int_\alpha\lambda_\Omega(w)|dw|.
$$
One of the most important properties of the hyperbolic metric is the
{\it principle of hyperbolic metric}, which asserts the monotonicity
$\lambda_\Omega(w)\ge\lambda_{\Omega_0}(w)$ and thus
$h_\Omega(w_1,w_2)\ge h_{\Omega_0}(w_1,w_2)$
for $w, w_1, w_2\in\Omega\subset\Omega_0$
({\it cf.}~\cite[III.3.6]{Nev:Ein}).
For basic facts about hyperbolic metrics, we refer to recent textbooks
\cite{KL:hg} by Keen and Lakic or a survey paper \cite{BM07} by Beardon and Minda
as well as classical book \cite{Ahlfors:conf} by Ahlfors.
For instance,
\begin{equation}\label{eq:hdist}
h_\uhp(\zeta_1, \zeta_2)=\arctanh\left|\frac{\zeta_1-\zeta_2}{\zeta_1-\overline{\zeta_2}}\right|,
\quad \zeta_1, \zeta_2\in\uhp,
\end{equation}
where $\arctanh x=\frac12\log\frac{1+x}{1-x}$ for $0\le x<1$
(see \cite[Chap.~7]{Beardon:disc} for instance).
Therefore, $h_\Omega(w_1,w_2)$ can be expressed via the universal covering
projection $w=p(\zeta)$ as follows (see \cite[Theorem 7.1.3]{KL:hg}):
$$
h_\Omega(w_1,w_2)=\min_{\zeta_2\in p\inv(w_2)} h_\uhp(\zeta_1,\zeta_2)
=\min_{\gamma\in\Gamma} h_\uhp(\zeta_1,\gamma(\zeta_2)),
$$
where $\zeta_1\in p\inv(w_1), \zeta_2\in p\inv(w_2)$ and $\Gamma=\{\gamma\in
\PSL(2,\R): p\circ\gamma=p\}$ is the covering transformation group
of $p:\uhp\to\Omega$ (also called a Fuchsian model of $\Omega$).

It is, however, difficult to obtain an explicit expression
of $\lambda_\Omega(w)$ or $h_\Omega(w_1,w_2)$ for a general hyperbolic domain $\Omega$
because a concrete form of its universal covering projection is not known
except for several special domains.
(It is not easy even for a simply connected domain because it is hard to
find its Riemann mapping function in general.)
Therefore, as a second choice, estimates of $\lambda_\Omega(w)$ are useful.
Indeed, Beardon and Pommerenke \cite{BP78} supplied a general but concrete bound
for $\lambda_\Omega(w).$
However, it is still difficult to estimate the induced hyperbolic distance $h_\Omega(w_1,w_2)$
due to complexity of the fundamental group of $\Omega.$

On the other hand, an explicit bound for the hyperbolic distance may be of importance.
For instance, if $f:\Omega\to X$ is a holomorphic map between hyperbolic domains,
then the principle of hyperbolic metric yields the inequality
\begin{equation}\label{eq:pr1}
h_X(f(z),f(z_0))\le h_\Omega(z,z_0),\quad z_0,z\in\Omega,
%=\arth\left|\frac{z-z_0}{1-\overline{z_0}z}\right|,
\end{equation}
which contains some important information about the function $f.$
As a maximal hyperbolic plane domain, the thrice-punctured sphere
(the twice-punctured plane) $\C_{0,1}:=\C\setminus\{0,1\}$ is particularly important.
Letting $X=\C_{0,1},$ the inequality \eqref{eq:pr1}
leads to Schottky's theorem when $\Omega$ is the unit disk ({\it cf.}~\cite{Hempel79}, \cite{Hempel80}),
and it leads to the big Picard theorem
when $\Omega$ is a punctured disk ({\it cf.}~\cite[\S 1-9]{Ahlfors:conf}).
Though the hyperbolic density of $\C_{0,1}$ was essentially computed by
Agard \cite{Agard68} (see also \cite{SV05}) and the holomorphic universal
covering projection of $\uhp$ onto $\C_{0,1}$ is known as an elliptic modular function
(see Section 2 below and, e.g., \cite[p.~279]{Ahlfors:ca} or \cite[Chap.~VI]{Nehari:conf}),
we do not have any convenient expression of the hyperbolic distance $h_{\C_{0,1}}(w_1,w_2)$ except for
special configurations of the points $w_1,w_2$ (see, for instance, \cite[Lemma 3.10]{SV01},
\cite[Lemma 5.1]{SV05}).

%More than ten years ago, Matti Vuorinen asked the first author to find a concrete but equivalent
%distance, or even a bound for $h_{\C_{0,1}}(w_1,w_2),$ where $\C_{0,1}=\C\setminus\{0,1\}.$
In this paper, we consider punctured spheres $X=\sphere\setminus\{a_1,\dots,a_n\}.$
This is still general enough in the sense that for any hyperbolic domain
$\Omega\subset\sphere,$ there exists a sequence of punctured spheres $X_k\supset\Omega~(k=1,2,\dots)$
such that $\lambda_{X_k}\to \lambda_\Omega$ locally uniformly on $\Omega$ as $k\to\infty$
(see \cite[\S 5]{BR86}).
We also note that Rickman \cite{Rickman84} constructed a conformal metric on a
punctured sphere of higher dimensions to show a Picard-Schottky type result for
quasiregular mappings.
Our main purpose of this paper is to give a distance function $d_X(w_1,w_2)$
on the punctured sphere $X$ which can be computed (or estimated) more easily than the hyperbolic
distance $h_X(w_1,w_2)$ but still comparable with it by concrete bounds.
To this end, we first propose a distance function $D(z_1,z_2)$
on $0<|z|\le e\inv$ given by the formula
\begin{equation} \label{eq:D}
D(z_1, z_2)=\frac{2\sin(\theta/2)}{\max\{\log (1/|z_1|),\, \log (1/|z_2|)\}}
+\left|\log\log\frac1{|z_2|}-\log\log\frac1{|z_1|}\right|,
\end{equation}
where $\theta=|\arg (z_2/z_1)| \in [0,\pi]$.
In Section 3, we will show that $D(z_1,z_2)$ is indeed a distance function
and comparable with $h_{\D^*}(z_1,z_2)$ on the set $0<|z|\le e\inv.$
Note also that $D(z_1, z_2)=e|z_1-z_2|$ for $|z_1|=|z_2|=e\inv.$

As another extremal case, the punctured disk $\D^*$ is also important.
Here and hereafter, we set
$\D(a,r)=\{z\in\C: |z-a|<r\}, \bD(a,r)=\{z\in\C: |z-a|\le r\},
~\D^*(a,r)=\D(a,r)\setminus\{a\}$
for $a\in\C$ and $0<r<+\infty$ and
%$\D_r=\D(0,r),~ \D_r^*=\D^*(0,r),~
$\D=\D(0,1)$ and $\D^*=\D^*(0,1).$
In this case, a holomorphic universal covering projection $p:\uhp\to\D^*$ is
given by $z=p(\zeta)=e^{\pi i\zeta}$ and the hyperbolic density is expressed by
\begin{equation*}%\label{hyper density in D star}
\lam_{\mathbb{D}^*}(z)=\frac{1}{2|z| \log (1/|z|)}.
\end{equation*}
A concrete formula of $h_{\D^*}(z_1,z_2)$ can also be given but its form
is not so convenient ({\it cf.} \eqref{eq:hd} below).

In order to understand the hyperbolic distance $h_X(w_1,w_2)$ when
one of $w_1,w_2$ is close to a puncture,
we should take a careful look at the hyperbolic geodesic nearby the puncture.
In Section 2, we investigate it by making use of an elliptic modular function as
well as the punctured unit disk model.
Section 3 is devoted to the study of the function $D(z_1,z_2).$
In particular, we show that $D$ gives a distance on $0<|z|\le e\inv$
and compare with the hyperbolic distances $h_{\D^*}(z_1,z_2)$ of $\D^*$
and $h_{\C_{0,1}}(z_1,z_2)$ of $\C_{0,1}.$
As an application of the function $D(z_1,z_2),$ in Section 4,
we will construct a distance function $d_X(w_1,w_2)$ on $n$-times
punctured spheres $X$ which are comparable with the hyperbolic distance $h_X(w_1,w_2).$
We will summarise our main results in Theorem \ref{thm:main}.
Unfortunately, $d_X(w_1,w_2)$ is not very easy to compute because we have to
take an infimum in the definition.
In Section 5, we introduce yet another quantity $e_X(w_1,w_2),$
which can be computed without taking an infimum though it is no longer a distance function on $X.$
We will show our main result that $d_X(w_1,w_2)$ and $e_X(w_1,w_2)$ are both comparable
with the hyperbolic distance $h_X(w_1,w_2)$ in a quantitative way in Section 5.

We would finally thank Matti Vuorinen for posing, more than ten years ago, the problem of finding a
quantity comparable with the hyperbolic distance on $\C_{0,1}$ and for helpful suggestions.

\section{Hyperbolic geodesics near the puncture}

In order to estimate the hyperbolic distance $h_X(w_1,w_2)$
of a punctured sphere $X=\sphere\setminus\{a_1,\dots,a_n\},$
we have to investigate the behaviour of a hyperbolic geodesic
joining two points near a puncture.
Here and in what follows, a curve $\alpha$ joining $w_1$ and $w_2$
in a hyperbolic domain $\Omega$
is called a {\it hyperbolic geodesic} if $\ell_\Omega(\alpha)\le\ell_\Omega(\beta)$
whenever $\beta$ is a curve joining $w_1$ and $w_2$
which is homotopic to $\alpha$ in $\Omega.$
In particular, $\alpha$ is called {\it shortest} if $\ell_\Omega(\alpha)
=h_\Omega(w_1,w_2).$
Note that the shortest hyperbolic geodesic is not unique in general.
Our basic model for that is the punctured disk $\D^*.$
In this case, we have precise information about the hyperbolic geodesic.

\begin{lem}\label{lem:geo}
Let $z_1, z_2\in\D^*$ with $\eta_1=-(\log|z_1|)/\pi \ge\eta_2=-(\log|z_2|)/\pi.$
Then a shortest hyperbolic geodesic $\beta$ joining $z_1$ and $z_2$ in $\D^*$
is contained in the set $e^{-\pi\delta}|z_1|\le |z|\le |z_2|,$ where
$$
\delta=\sqrt{\eta_1^2+\frac1{4}}-\eta_1
=\frac1{4\big(\eta_1+\sqrt{\eta_1^2+1/4}\big)}.
$$
\end{lem}

\begin{pf}
First note that the function $p(\zeta)=e^{\pi i\zeta}$ is a universal covering
projection of the upper half-plane $\uhp$ onto $\D^*$ with period $2.$
We may assume that $z_1\in(0,1)$ and $\theta=(\arg z_2)/\pi\in(0,1].$
Then $p(i\eta_1)=z_1,~ p(i\eta_2+\theta)=z_2,$ and
$h_\uhp(i\eta_1,i\eta_2+\theta)=h_{\D^*}(z_1,z_2).$
Let $\tilde\beta$ be the hyperbolic geodesic joining $\zeta_1=i\eta_1$
and $\zeta_2=i\eta_2+\theta$ in $\uhp.$
Recall that $\tilde\beta$ is part of the circle orthogonal to the real axis.
If we fix $\eta_1,$ the possible largest imaginary part of $\tilde\beta$ is attained
when $\eta_2=\eta_1$ and $\theta=1.$
Therefore, $\Im \zeta\le \eta_1+\delta$ for $\zeta\in\tilde\beta,$ where
$\delta=\sqrt{\eta_1^2+1/4}-\eta_1.$
Hence, we conclude that $\beta=p(\tilde\beta)$ is contained in the closed annulus
$e^{-\pi(\eta_1+\delta)}=e^{-\pi\delta}|z_1|\le |z|\le |z_2|$ as required.
\end{pf}

By the proof, we observe that the above constant $\delta$ is sharp.
Note here that $\delta$ is decreasing in $\eta_1$ and that $0<\delta<\frac1{8\eta_1}.$

In the above theorem, we see that the subdomain $\D^*(0,\rho)$ of $\D^*$
with $0<\rho<1$ is hyperbolically convex.
This is also true in general.
Indeed, the following result is a special case of Minda's reflection
principle \cite{Minda87}
(apply his Theorem 6 to the case when $\overline R=\overline\Delta$).

\begin{lem}\label{lem:Minda}
Let $\Omega$ be a hyperbolic subdomain of $\C$ and let $\Delta$ be
an open disk centered at a point $a\in\C\setminus\Omega.$
Suppose that $I(\Omega\setminus\Delta)\subset\Omega,$ where $I$ denotes the
reflection in the circle $\partial\Delta.$
Then $\Delta\cap\Omega$ is hyperbolically convex in $\Omega.$
\end{lem}

In particular, we have

\begin{cor}\label{cor:Minda}
Let $\Delta$ be an open disk centered at a puncture
$a$ of a hyperbolic punctured sphere $X\subset\C$
with $\Delta^*=\Delta\setminus\{a\}\subset X.$
Then $\Delta^*$ is hyperbolically convex in $X.$
\end{cor}

As another extremal case, we now consider the thrice-punctured sphere
$\C_{0,1}.$
It is well known that the elliptic modular function, which is denoted by
$J(\zeta),$ on the upper half-plane $\uhp=\{\zeta:\Im\zeta>0\}$ serves as
a holomorphic universal covering projection onto $\C_{0,1}.$
The reader can consult \cite[Chap.~VI]{Nehari:conf} for general facts about the function
$J$ and related functions.
The covering transformation group is the modular group $\Gamma(2)$ of level 2;
namely, $\Gamma(2)=\{A\in\SL(2,\Z): A\equiv I~\mod 2\}/\{\pm I\}.$
Since $J$ is periodic with period 2, it factors into $J=Q\circ p,$
where $p(\zeta)=e^{\pi i\zeta}$ as before and $Q$ is an intermediate covering projection
of $\D^*$ onto $\C_{0,1}.$
Since $J(\zeta)\to0$ as $\eta=\Im\zeta\to+\infty,$ the origin is a removable
singularity of $Q(z)$ and indeed the function $Q(z)$ has the following representations
(see \cite[VI.6]{Nehari:conf} or \cite[Theorem 14.2.2]{KL:hg}):
\begin{align*}
Q(z)&=16z\prod_{n=1}^\infty\left(\frac{1+z^{2n}}{1+z^{2n-1}}\right)^8
=16z\left[\frac{\sum_{n=0}^\infty z^{n(n+1)}}{1+2\sum_{n=1}^\infty z^{n^2}}\right]^4
\\
&=16(z-8z^2+44z^3-192z^4+718z^5-\cdots).
\end{align*}
We also remark that the function $Q(z)$ has been recently used to improve
coefficient estimates of univalent harmonic mappings on the unit disk
in Abu Muhanna, Ali and Ponnusamy \cite{AAP17}.
By its form, $Q(z)$ is locally univalent at $z=0.$
In fact, we can show the following.

\begin{lem}\label{lem:univ}
The function $Q(z)$ is univalent on the disk $|z|<e^{-\pi/2}\approx0.20788.$
The radius $e^{-\pi/2}$ is best possible.
\end{lem}

\begin{pf}
Suppose that $Q(z_1)=Q(z_2)$ for a pair of points $z_1,z_2$ in the disk $|z|<e^{-\pi/2}.$
Take the unique point $\zeta_l=\xi_l+i\eta_l\in p\inv(z_l)\subset\uhp$ such that $-1<\xi_l\le 1$
and $\eta_l>1/2$ for $l=1,2.$
We now recall the well-known fact that $\omega=\{\zeta\in\uhp: -1<\Re\zeta\le1,
|\zeta+1/2|<1/2, |\zeta-1/2|\le1/2\}$ is a fundamental set for the modular group
$\Gamma(2)$ of level 2 (see \cite[Chap.7, \S 3.4]{Ahlfors:ca}).
In other words, $J(\omega)=\C_{0,1}$ and no pairs of distinct points in $\omega$
have the same image under the mapping $J.$
We now note that $\zeta_l\in \omega$ for $l=1.2.$
Since $J(\zeta_1)=Q(z_1)=Q(z_2)=J(\zeta_2),$ we conclude that $\zeta_1=\zeta_2.$
Hence, $z_1=z_2,$ which implies that $Q(z)$ is univalent on $|z|<e^{-\pi/2}.$
To see sharpness, we consider
the pair of points $\zeta_1=(1+i)/2$ and $\zeta_2=(-1+i)/2.$
Since $\zeta_2=T(\zeta_1)$ for the modular transformation $T(\zeta)=\zeta/(-2\zeta+1)$ in $\Gamma(2),$
we have $J(\zeta_1)=J(\zeta_2).$
Thus $Q(ie^{-\pi/2})=Q(-ie^{-\pi/2}).$
\end{pf}

We remark that the formula $Q(ie^{-\pi/2})=2$ is valid.
Indeed, by recalling the functional identity $Q(-z)=Q(z)/(Q(z)-1)$
(see \cite[(92) in p.~328]{Nehari:conf}),
we have $Q(ie^{-\pi/2})=Q(-ie^{-\pi/2})
=Q(ie^{-\pi/2})/(Q(ie^{-\pi/2})-1)$, which implies $Q(ie^{-\pi/2})=2.$

We now recall the growth theorem for a normalized univalent analytic function $f(z)=z+a_2z^2+\dots$
on $|z|<1$ (see \cite[\S 5-1]{Ahlfors:conf} for instance):
$$
\frac{r}{(1+r)^2}\le |f(z)|\le \frac{r}{(1-r)^2},\quad |z|<1.
$$
Applying this result to $f(z)=e^{\pi/2}Q(ze^{-\pi/2})/16,$ we obtain the following estimates:
\begin{equation}\label{eq:Q}
\frac{16r}{(1+re^{\pi/2})^2}\le |Q(z)|\le\frac{16r}{(1-re^{\pi/2})^2},
\quad r=|z|<e^{-\pi/2}.
\end{equation}
Observe that the lower bound in \eqref{eq:Q} tends to $4e^{-\pi/2}$
as $r\to e^{-\pi/2}.$
Hence $\D(0,4e^{-\pi/2})\subset Q(\D(0,e^{-\pi/2})).$
Note that $4e^{-\pi/2}\approx 0.8315.$

\begin{lem}\label{lem:N}
Let $\C_{0,1}.$
For $w_1, w_2\in \C_{0,1}$ with $|w_1|, |w_2|\le\rho\le 4e^{-\pi/2},$
a shortest hyperbolic geodesic $\alpha$ joining
$w_1,w_2$ in $\C_{0,1}$ is contained in the closed annulus
$$
\min\{|w_1|,|w_2|\}e^{-K}\le |w|\le \max\{|w_1|,|w_2|\},
$$
where $K=K(\rho)>0$ is a constant depending only on $\rho.$
\end{lem}

\begin{pf}
The right-hand inequality follows from Corollary \ref{cor:Minda}.
We now show the left-hand inequality.
Take $0<r\le e^{-\pi/2}$ so that $16r/(1+re^{\pi/2})^2=\rho$ and put
$\mu=re^{\pi/2}\le 1.$
Note that $r$ and $\mu$ can be computed by the formula
$$
\mu=re^{\pi/2}
=\rho\inv e^{-\pi/2}\big(8-\rho e^{\pi/2}-4\sqrt{4-\rho e^{\pi/2}}\big).
$$
By \eqref{eq:Q}, we can choose
$z_j\in \D^*$ with $|z_j|\le r$ and $Q(z_j)=w_j$ for $j=1,2$
in such a way that a lift $\beta$ of the curve $\alpha$ joins
$z_1$ and $z_2$ in $\D^*$ via the covering map $Q.$
We may assume that $|z_1|\le |z_2|.$
It follows from Lemma \ref{lem:geo} that $\beta$ is contained in
the annulus $|z_1|e^{-\pi\delta}\le |z|\le |z_2|,$ where $\delta$
is given in the lemma with $\eta_1=-(\log r)/\pi \le -(\log|z_1|)/\pi.$
Hence, by \eqref{eq:Q}, the curve $\alpha=Q(\beta)$ is contained
in the annulus
$$
\frac{16|z_1|e^{-\pi\delta}}{(1+|z_1|e^{\pi/2-\pi\delta})^2}\le
|w|\le\frac{16|z_2|}{(1-|z_2|e^{\pi/2})^2}.
$$
By \eqref{eq:Q} and $|z_1|e^{\pi/2}\le\mu,$ we see that
$$
|w_1|\le \frac{16|z_1|}{(1-|z_1|e^{\pi/2})^2}
=\frac{16|z_1|e^{-\pi\delta}}{(1+|z_1|e^{\pi/2-\pi\delta})^2} \cdot
\frac{e^{\pi\delta}(1+|z_1|e^{\pi/2-\pi\delta})^2}{(1-|z_1|e^{\pi/2})^2}
\le \frac{16|z_1|e^{-\pi\delta}}{(1+|z_1|e^{\pi/2-\pi\delta})^2} e^K,
$$
where
\begin{equation}\label{eq:K}
K=K(\rho)=\pi\delta+2\log\frac{1+\mu e^{-\pi\delta}}{1-\mu}.
\end{equation}
Thus we have seen that $\alpha$ is contained in the set $|w_1|e^{-K}\le |w|.$
\end{pf}

As an immediate consequence, we obtain the following result.

\begin{cor}\label{cor:rho}
Let $0<\sigma\le\rho\le 4e^{-\pi/2}.$
A shortest hyperbolic geodesic $\alpha$ joining
$z_1,z_2$ in $\C_{0,1}$ with $|z_1|,|z_2|\ge\sigma$ does not intersect the disk
$\Delta=\D(0,\sigma e^{-K}),$ where $K=K(\rho)>0$ is the constant in Lemma \ref{lem:N}.
\end{cor}

\begin{pf}
Suppose that $\alpha$ intersects the disk $\Delta.$
Then we can choose a subarc $\alpha_0$ of $\alpha$ such that
$\alpha_0$ intersects $\Delta$ and that
both endpoints $w_1,w_2$ of $\alpha_0$ have modulus $\sigma.$
Applying Lemma \ref{lem:N} to $\alpha_0$ yields a contradiction.
\end{pf}

When $\rho=\rho_0=e\inv,$ we compute $r_0=e^{1-\pi}\big(8-e^{\pi/2-1}-4\sqrt{4-e^{\pi/2-1}}\big)
\approx 0.0301441$ and $\mu_0=r_0 e^{\pi/2}\approx 0.145007.$
Also, we have $\delta_0=\sqrt{\eta_0^2+1/4}-\eta_0\approx 0.107007$ with
$\eta_0=-(\log r_0)/\pi\approx 1.11465.$
Then $K_0=K(e\inv)=\pi\delta_0+2\log[(1+\mu_0e^{-\pi\delta_0})/(1-\mu_0)]\approx 0.846666$
and $e^{K_0}\approx 2.33186.$
Thus we obtain the following statement as a special case of  Corollary \ref{cor:rho}.

\begin{cor}\label{cor:N}
Let $z_1, z_2$ be two points in $\C_{0,1}$ with $|z_1|, |z_2|\ge\sigma$ for
some number $\sigma\in(0, e\inv].$
Then, a shortest hyperbolic geodesic $\alpha$ joining
$z_1,z_2$ in $\C_{0,1}$ does not intersect the disk
$|z|<\sigma e^{-K_0},$ where $K_0=K(e\inv)\approx 0.846666.$
\end{cor}

\section{Basic properties of the distance function $D(w_1,w_2)$}

First we show the following result in the present section.
Let $E^*=\{z: 0<|z|\le e\inv\}.$

\begin{lem} \label{D distance theorem}
The function $D(z_1, z_2)$ given by \eqref{eq:D} is a distance function
on the set $E^*.$
\end{lem}

\begin{pf}
First we note that $D(z_1, z_2)$ can also be described by
$$
D(z_1,z_2)=\frac{|\zeta_1-\zeta_2|}{\max\{\tau_1, \tau_2\}}
+|\log\tau_1-\log\tau_2|,
$$
where $\tau_j=\log(1/|z_j|)$ and $\zeta_j=z_j/|z_j|$ for $j=1,2.$
It is easy to see that $D(z_1,z_2)=D(z_2,z_1)$
and that $D(z_1,z_2) \geq 0,$ where equality holds if and only if $z_1=z_2$.
It remains to verify the triangle inequality.
Our task is to show that the inequality
$\Delta:=D(z_1,z)+D(z,z_2)-D(z_1, z_2)\ge0$ for $z\in E^*.$
Set $\tau=\log(1/|z|)$ and $\zeta=z/|z|.$
We may assume that $\tau_2\le\tau_1.$
Then
\begin{equation*} %\label{integral D}
D(z_1, z_2)
=\frac{|\zeta_1-\zeta_2|}{\tau_1}
+\log\tau_1-\log\tau_2.
\end{equation*}

First we assume that $\tau\le\tau_1.$
Since $\tau_1\ge\max\{\tau, \tau_2\},$ we have
\begin{align*}
D(z_1,z_2)&\le\frac{|\zeta_1-\zeta|+|\zeta-\zeta_2|}{\tau_1}
+|\log\tau_1-\log\tau|+|\log\tau-\log\tau_2| \\
&= D(z_1,z)+D(z,z_2).
\end{align*}
Secondly, we assume that $\tau>\tau_1.$
Then, in a similar manner, we have
\begin{align*}
\Delta&\ge \left(\frac1\tau-\frac1{\tau_1}\right)
|\zeta_1-\zeta_2|+2\log\tau-2\log\tau_1 \\
&\ge \frac2\tau-\frac2{\tau_1}+2\log\tau-2\log\tau_1.
\end{align*}
Since the function $f(x)=2/x+2\log x$ is increasing in $1\le x<+\infty,$
we have $\Delta\ge0$ as required.
\end{pf}

Next we compare $D(z_1,z_2)$ with the hyperbolic distance $h_{\D^*}(z_1,z_2)$
of $\D^*$ on the set $E^*.$

\begin{thm} \label{thm:ct}
The distance function $D(z_1, z_2)$ given by \eqref{eq:D} satisfies
\begin{equation}\label{eq:cp}
\frac{4}{\pi}\,h_{\D^*}(z_1,z_2) \leq D(z_1,z_2) \leq
M_0\, h_{\C_{0,1}}(z_1,z_2)
\end{equation}
for $0<|z_1|, |z_2|\le e\inv.$
Here, $\C_{0,1}=\C\setminus\{0,1\},$
$M_0$ is a positive constant with $M_0<24$ and the constant $4/\pi$ is sharp.
\end{thm}

We remark that $4/\pi\approx 1.27324.$

\begin{pf}
We consider the quantity
$$
D'(z_1, z_2)=\frac{\theta}{\max\{\log (1/|z_1|),\, \log (1/|z_2|)\}}
+\left|\log\log\frac1{|z_2|}-\log\log\frac1{|z_1|}\right|
$$
for $z_1, z_2\in \D^*,$ where $\theta=|\arg(z_2/z_1)|\in[0,\pi].$
Since $2x/\pi\le\sin x\le x$ for $0\le x\le\pi/2,$ we can easily obtain
$$%\be \label{l and D}
\frac{2}{\pi} D'(z_1,z_2) \leq D(z_1,z_2) \leq D'(z_1,z_2)
$$
for $z_1, z_2 \in E^*$.
We show now the inequality
\begin{equation}\label{eq:L}
2h_{\D^*}(z_1,z_2) \leq D'(z_1,z_2)
\end{equation}
for $z_1, z_2\in E^*.$
Combining these two inequalities, we obtain the first inequality
in \eqref{eq:cp}.

Without loss of generality, we may assume that
$\arg z_1 =0$, $\arg z_2 = \theta \in [0, \pi]$,
and $\tau_2 \le\tau_1,$ where $\tau_j=\log(1/|z_j|).$
Recall that $p(\zeta)=e^{\pi i\zeta}$ is a holomorphic universal covering
projection of the upper half-plane $\uhp$ onto $\D^*.$
Let $\zeta_1=i\tau_1/\pi$ and $\zeta_2=(\theta+i\tau_2)/\pi.$
Then $p(\zeta_j)=z_j$ for $j=1,2$ and, by \eqref{eq:hdist},
\begin{equation}\label{eq:hd}
h_{\D^*}(z_1,z_2)=h_\uhp(\zeta_1,\zeta_2)
=\arctanh\left|\frac{\zeta_1-\zeta_2}{\zeta_1-\overline{\zeta_2}}\right|
=\arctanh\sqrt{\frac{\theta^2+(\tau_1-\tau_2)^2}{\theta^2+(\tau_1+\tau_2)^2}}.
%=\frac12\log \frac{|z_1-\overline{z_2}|+|z_1-z_2|}%
%{|z_1-\overline{z_2}|-|z_1-z_2|} \\
%&=\frac12\log \frac{\sqrt{\theta^2+(\tau_1+\tau_2)^2}
%+\sqrt{\theta^2+(\tau_1-\tau_2)^2}}{\sqrt{\theta^2+(\tau_1+\tau_2)^2}
%-\sqrt{\theta^2+(\tau_1-\tau_2)^2}}.
\end{equation}
We consider the function
\begin{align*}
G(\theta, \tau_1,\tau_2)&= 2h_{\D^*}(z_1,z_2)-D'(z_1,z_2) \\
&=2\arctanh\sqrt{\frac{\theta^2+(\tau_1-\tau_2)^2}{\theta^2+(\tau_1+\tau_2)^2}}
-\frac\theta{\tau_1}-\log\tau_1+\log\tau_2
\end{align*}
on $0\le\theta\le\pi, 1\le\tau_2\le\tau_1.$
A straightforward computation yields
$$
\frac{\partial G}{\partial \tau_2}=\frac{1}{\tau_2}
-\frac{\tau_1^2-\tau_2^2+\theta^2}%
{\tau_2\sqrt{(\theta^2+\tau_2^2+\tau_1^2)^2-4\tau_1^2 \tau_2^2}}.
$$
Since
$$
\left[(\theta^2+\tau_2^2+\tau_1^2)^2-4\tau_1^2 \tau_2^2\right]
-(\tau_1^2-\tau_2^2+\theta^2)^2
=4\tau_2^2\theta^2>0
$$
for $0<\theta<\pi,$
we see that $\partial G/\partial\tau_2>0$ and therefore
$G(\theta,\tau_1,\tau_2)$ is increasing in $1\le \tau_2\le\tau_1$
so that
$$
G(\theta,\tau_1,1)<G(\theta,\tau_1,\tau_2)<G(\theta,\tau_1,\tau_1)
$$
for $1<\tau_2<\tau_1.$
We observe that $G(\theta,\tau_1,\tau_1)=2f(2\tau_1/\theta),$ where
$$
f(x)=\arctanh\frac1{\sqrt{1+x^2}}-\frac1x
=\frac{g(x)-1}{x},
$$
and $g(x)=x\,\arctanh(1/\sqrt{1+x^2}).$
%Note that the range of the quantity $x=2\tau_1/\theta$ for
%$1<\tau_1, 0<\theta<\pi$ is the interval $(2/\pi,+\infty).$
Since $g'(x)=\arctanh(1/\sqrt{1+x^2})-1/\sqrt{1+x^2}>0$
for $x>0,$
we have $g(x)<g(+\infty)=1.$
Hence $G(\theta,\tau_1,\tau_1)<0$ for $0<\theta<\pi, \tau_1>1.$
We have thus proved the inequality $2h_{\D^*}(z_1,z_2)\le D'(z_1,z_2).$
Since $h_{\D^*}(e^{-\tau},-e^{-\tau})=h_\uhp(i\tau,i\tau+\pi)
=\arth(\pi/\sqrt{\pi^2+4\tau^2}),$
$$
\frac{h_{\D^*}(e^{-\tau},-e^{-\tau})}{D(e^{-\tau},-e^{-\tau})}
=\frac\tau 2\,\arctanh\frac{\pi}{\sqrt{\pi^2+4\tau^2}}\to\frac{\pi}4
$$
as $\tau\to+\infty.$
Hence the constant $4/\pi$ is sharp in \eqref{eq:cp}.

Finally, we show the second inequality in \eqref{eq:cp}.
Assume that $0<|z_1|\le|z_2|\le e\inv, \theta_0=\arg(z_2/z_1)\in[0,\pi]$
and $\tau_j=-\log|z_j|\ge 1$ for $j=1,2.$
We now show the following two inequalities to complete the proof:
\begin{equation}\label{eq:two}
\frac{\theta_0}{\tau_1}\le M_1H
\aand
\log\frac{\tau_1}{\tau_2}\le M_2H
\end{equation}
for some constants $M_1$ and $M_2,$ where $H=h_{\C_{0,1}}(z_1,z_2).$
Let $\alpha$ be a shortest hyperbolic geodesic joining $z_1$ and $z_2$ in $\C_{0,1}.$
Then, by Corollary \ref{cor:N}, $\alpha$ is contained in the annulus
$|z_1|e^{-K_0}\le|z|\le |z_2|.$
Thus $-\log|z|\le K_0+\tau_1$ for $z\in\alpha.$
We recall the following lower estimate of $\lambda(z)=\lambda_{\C_{0,1}}(z)$
(see \cite{Hempel79}):
\begin{equation}\label{eq:Hempel}
\frac1{2|z|(C_0+|\log|z||)}\le \lambda(z),\quad z\in \C_{0,1},
\end{equation}
where $C_0=1/2\lambda(-1)=\Gamma(1/4)^4/4\pi^2\approx 4.37688.$
We remark that the bound in \eqref{eq:Hempel} is monotone decreasing in $|z|.$
Noting the inequality $|dz|\ge rd\theta$ for $z=re^{i\theta},$ we have
$$
H=\int_\alpha\lambda(z)|dz|
\ge \int_\alpha \frac{rd\theta}{2r(C_0-\log r)}
\ge \int_\alpha \frac{d\theta}{2(C_0+\tau_1+K_0)}
=\frac{\theta_0}{2(C_0+K_0+\tau_1)}.
$$
Since $C_0+K_0+\tau_1\le (C_0+K_0+1)\tau_1,$ we obtain the
first inequality in \eqref{eq:two} with $M_1=2(C_0+K_0+1).$
Similarly, by using $|dz|\ge |dr|$ for $z=re^{i\theta}=e^{-t+i\theta},$
we obtain
$$
H\ge \int_\alpha\frac{|dr|}{2r(C_0-\log r)}
\ge\int_{\tau_2}^{\tau_1}\frac{dt}{2(C_0+t)}
\ge \frac1{2(C_0+1)}\log\frac{\tau_1}{\tau_2}.
$$
Hence we have the second inequality in \eqref{eq:two} with
$M_2=2(C_0+1).$
Combining the two inequalities in \eqref{eq:two}, we get
$$
D(z_1,z_2)\le D'(z_1,z_2)\le (M_1+M_2)H.
$$
Hence, $M_0=M_1+M_2=4C_0+4+2K_0\approx 23.2008$ works.
\end{pf}

\section{Construction of a distance function on $n$-times punctured sphere}

In this section, we construct a distance function on an $n$-times
punctured sphere $X=\sphere\setminus\{a_1,\dots, a_n\}.$
As we noted, $X$ is hyperbolic if and only if $n\ge3.$
Thus we will assume that $n\ge3$ in the sequel.
After a suitable M\"{o}bius transformation, without much loss of generality,
we may assume that $0,1,\infty\in\sphere\setminus X$
and $a_n =\infty$.
Let
$$
\tilde\rho_j=\begin{cases}
\displaystyle
\min_{1\le k<n, k\ne j}|a_k-a_j|&\quad \text{for}~j=1,2,\dots, n-1, \\
\displaystyle
\max_{1<k<n}|a_k| &\quad \text{for}~j=n
\end{cases}
$$
and $\rho_j=\tilde\rho_j/e$ for $j=1,2,\dots, n-1$ and $\rho_n=e\tilde\rho_n.$
Since $a_1=0,$ we have $e\rho_j\le|a_j|$ for $1<j<n$ and
$e\rho_j\le\max |a_k|=\rho_n/e$ for $j<n.$
Set $E_{j}=\bD(a_j,\rho_j), E^*_j=E_j\setminus\{a_j\}$ and
$\Delta_{j}=\D^*(a_j,\tilde\rho_j)$
for $j=1,\dots, n-1,$ and set
$E_{n}=\{w\in\sphere: |w|\ge \rho_n\},~E^*_n=E_n\setminus\{\infty\}$ and
$\Delta_{n}=\{w\in\C: |w|>\tilde\rho_n\}.$
It is easy to see that the Euclidean distance between $E_j$'s are computed and estimated by
$$
\dist(E_j,E_k)=|a_j-a_k|-\rho_j-\rho_k
\ge(e-2)\max\{\rho_j,\rho_k\}
$$
for $j,k<n,$ and
$$
\dist(E_j,E_n)=\rho_n-|a_j|-\rho_j\ge
(1-e\inv-e^{-2})\rho_n\ge (e^2-e-1)\rho_j>(e-2)\rho_j
$$
for $1\le j<n.$
In particular, $E^*_j$'s are mutually disjoint.
Noting the inequality $1-e\inv-e^{-2}\approx0.49678<e-2,$
we also have the estimate
\begin{equation}\label{eq:dist}
\dist(E_j,E_k)\ge (1-e\inv-e^{-2})\rho_j
\end{equation}
for any pair of distinct $j,k.$
Note also that $\Delta_j\subset X$ for $j=1,\dots, n.$
Finally, let $W=X\setminus (E^*_1\cup\dots\cup E^*_n)$.

We are now ready to construct a distance function on $X.$
Set
\begin{equation*}
D_j(w_1, w_2)=
\begin{cases}
\rho_jD((w_1-a_j)/\tilde\rho_j,(w_2-a_j)/\tilde\rho_j)
& \text{if}~j=1,\dots, n-1, \\
\rho_n D(\tilde\rho_n/w_1,\tilde\rho_n/w_2) & \text{if}~j=n
\end{cases}
\end{equation*}
for $w_1, w_2\in E^*_j,$ where $D(z_1,z_2)$ is given in \eqref{eq:D}.
By definition, we have $D_j(w_1,w_2)=|w_1-w_2|$ for $w_1$, $w_2 \in \pa E_j$.
By Lemma \ref{D distance theorem},
we know that $D_j(w_1, w_2)$ is a distance function on $E^*_{j}.$
We further define
\begin{equation} \label{eq:cm}
d_X(w_1, w_2)=
\begin{cases}
\displaystyle D_j(w_1, w_2) & \mbox{if\ }\ w_1,\, w_2 \in E^*_{j}, \\
\displaystyle \inf_{\zeta \in \pa E_j} (D_j(w_1, \zeta) +|\zeta-w_2|)
&\mbox{if\ }\ w_1 \in E^*_{j}, \,  w_2 \in W,  \\
\displaystyle \inf_{\zeta \in \pa E_j} (|w_1-\zeta| +D_j(\zeta,w_2))
&\mbox{if\ }\ w_1 \in W,\, w_2\in E^*_{j},  \\
\displaystyle \inf_{\substack{\zeta_1 \in \pa E_j\\ \zeta_2 \in \pa E_k}}
(D_j(w_1,\zeta_1) +|\zeta_1-\zeta_2|+D_k(\zeta_2, w_2) )
&\mbox{if\ }\ w_1 \in E^*_{j}, \,  w_2 \in E^*_{k},\, j\ne k,\\
\displaystyle |w_1-w_2| & \mbox{if\ }\ w_1,\, w_2 \in W
\end{cases}
\end{equation}
for $w_1, w_2\in X.$
Note that the infima in the above definition can be replaced by minima.
Then we have the following result.

\begin{lem}\label{lem:dist}
$d_X$ is a distance function on $X.$
\end{lem}

\begin{pf}
It is easy to see that $d_X(w_1,w_2)=d_X(w_2,w_1)$, and $d_X(w_1,w_2) \geq 0,$
where equality holds if and only if $w_1=w_2$, for $w_1, w_2 \in X$.
It remains to verify the triangle inequality:
$d_X(w_1,w_2)\le d_X(w_1,w_3)+d_X(w_3,w_2).$
According to the location of these points, we need to consider several cases.
For instance, we consider the case when $w_1 \in E^*_{j}$, $w_2 \in E^*_{k}$
and $w_3\in E^*_l$ for distinct $j,k,l.$
Then,
\ben
&&d_X(w_1,w_3)+d_X(w_3,w_2)\\
&=&\inf_{\substack{\zeta_1\in\partial E_j,\zeta_2 \in\partial E_k\\ \zeta_3,\zeta_4\in\pa E_l}}
(D_j(w_1,\zeta_1)+|\zeta_1-\zeta_3|+D_l(\zeta_3,w_3)+D_l(w_3,\zeta_4)
+|\zeta_4-\zeta_2|+D_k(\zeta_2, w_2))\\
&\geq& \inf(D_j(w_1,\zeta_1)+|\zeta_1-\zeta_3|+D_l(\zeta_3,\zeta_4)
+|\zeta_4-\zeta_2|+D_k(\zeta_2,w_2))\\
&=& \inf(D_j(w_1,\zeta_1)+|\zeta_1-\zeta_3|+|\zeta_3-\zeta_4|
+|\zeta_4-\zeta_2|+D_k(\zeta_2,w_2))\\
&\geq& \inf_{\zeta_1 \in\partial E_j,\zeta_2 \in\partial E_k}
(D_j(w_1,\zeta_1)+|\zeta_1-\zeta_2|+D_k(\zeta_2,w_2))= d_X(w_1,w_2).
\een
The other cases can be handled similarly and therefore will be omitted.
\end{pf}

We remark that we can construct a similar distance when $n=2.$
Let $a_1=0$ and $a_2=\infty$ and consider $X=\sphere\setminus\{a_1,a_2\}
=\C^*.$
Then, we set
$$
d_{\C^*}(w_1,w_2)=
\begin{cases}
D(w_1/e, w_2/e) &\text{if}~ 0<|w_1|\le 1, 0<|w_2|\le 1, \\
D(1/(ew_1), 1/(ew_2)) & \text{if}~ 1\le|w_1|, 1\le |w_2|, \\
\displaystyle
\inf_{|\zeta|=1} (D(w_1/e, \zeta/e)+D(1/(e\zeta),1/(ew_2)) &\text{if}~
 0<|w_1|\le 1, 1\le |w_2|, \\
\displaystyle
\inf_{|\zeta|=1} (D(1/(ew_1), 1/(e\zeta))+D(\zeta/e,w_2/e) &\text{if}~
 1\le |w_1|, 0<|w_2|\le 1,
\end{cases}
$$
where $D$ is defined by \eqref{eq:D}.
Then we can see that $d_{\C^*}$ is a distance function on $\C^*.$
The asymptotic behaviour of $d_{\C^*}$ near the punctures are rather
different from that of the quasi-hyperbolic distance $q$ on $\C^*$
since $q(w_1,w_2)=|\log|w_1|-\log|w_2||+O(1)$ as $w_1,w_2\to0$
(see \cite{MO86} for instance).

Our main result in the present paper is the following.
In the next section, we will prove it in a stronger form (Theorem \ref{thm:e}).

\begin{thm} \label{thm:main}
The distance function $d_X(w_1, w_2)$ given in \eqref{eq:cm}
on the $n$-times punctured sphere $X=\sphere\setminus\{a_1,\dots,a_n\}
\subset\C_{0,1}$ is comparable with the hyperbolic distance $h_X(w_1,w_2)$ on $X.$
\end{thm}

\section{Proof of Theorem \ref{thm:main}}
We recall that $X=\sphere\setminus \{a_1,\dots, a_n\}\subset\C_{0,1}$ with $a_n=\infty.$
The function $d_X$ defined in the previous section has the merit that it gives a distance on $X.$
On the other hand, it is not easy to compute the exact value of $d_X(w_1,w_2)$
for a given pair of points $w_1, w_2\in X.$
The following quantity can be a good substitute of $d_X(w_1,w_2)$ because
it is computed easily, though it is not necessarily a distance function:

\begin{equation*} \label{eq:p}
e_X(w_1, w_2)=
\begin{cases}
\displaystyle D_j(w_1, w_2) & \mbox{if\ }\ w_1,\, w_2 \in E^*_{j}, \\
\displaystyle  D_j(w_1, \zeta) +|\zeta-w_2| &\mbox{if\ }\ w_1 \in E^*_{j}, \,  w_2 \in W  \\
\displaystyle |w_1-\zeta| +D_j(\zeta,w_2) &\mbox{if\ }\ w_1 \in W,\, w_2\in E^*_{j},  \\
\displaystyle D_j(w_1,\zeta_1) +|\zeta_1-\zeta_2|+D_k(\zeta_2, w_2)
&\mbox{if\ }\ w_1 \in E^*_{j}, \,  w_2 \in E^*_{k},\, j\ne k,\\
\displaystyle |w_1-w_2| & \mbox{if\ }\ w_1,\, w_2 \in W
\end{cases}
\end{equation*}
for $w_1, w_2\in X,$
where $\zeta$ is the intersection point of the line segment $[w_1,w_2]$
with the circle $\partial E_j$ in the second case,
and $\zeta_1$ and $\zeta_2$ are the intersection points of the
line segments $[w_1,w_2]$ with the circles $\partial E_j$
and $\partial E_k,$ respectively, in the third case.
By definition, the inequality $d_X(w_1,w_2)\le e_X(w_1,w_2)$ holds obviously.
Theorem \ref{thm:main} now follows from the next result.

\begin{thm}\label{thm:e}
There exist positive constants $N_1$ and $N_2$ such that the following inequalities hold:
$$
N_1h_X(w_1,w_2)\le d_X(w_1,w_2)\le e_X(w_1,w_2)\le N_2h_X(w_1,w_2)
$$
for $w_1,w_2\in X=\sphere\setminus\{a_1,a_2,\dots,a_n\}\subset\C_{0,1}.$
\end{thm}

\begin{pf}
Assume that $a_1=0$ and $a_n=\infty$ as before.
We recall that $W=X\setminus(E^*_1\cup\dots\cup E^*_n).$
Since $\lambda_X(w)\delta_X(w)\le 1,$ we obtain
\begin{equation}\label{eq:m}
\lambda_X(w)\le \frac1m,\quad w\in \overline W,
\end{equation}
where $m=\min_{1\le j<n}\rho_j.$
We show the first inequality.
Fix $w_1, w_2\in X$ and
assume that $w_1\in E^*_j$ and $w_2\in W$ for some $j.$
We can deal with the other cases similarly and thus we will omit it.
We further assume, for a moment, that $j\ne n.$
By definition,
$$
d_X(w_1,w_2)=D_j(w_1,\zeta_0)+|\zeta_0-w_2|
$$
for some $\zeta_0\in\partial E_j.$
If the line segment $L=[\zeta_0,w_2]$ intersects $E^*_k$ for some
$k\ne j,$ we replace the part $L\cap\D(a_k,\rho_k)$ of $L$
by the shorter component of $\partial E_k\setminus L$
for each such $k.$
The resulting curve will be denoted by $L'.$
It is obvious from construction that the Euclidean length of
$L'$ is bounded by $\pi|\zeta_0-w_2|/2.$
Therefore, by \eqref{eq:m},
$$
|\zeta_0-w_2|\ge \frac2\pi \int_{L'}|dw|
\ge \frac{2m}{\pi}\int_{L'}\lambda_X(w)|dw|
\ge \frac{2m}{\pi}h_X(\zeta_0,w_2).
$$
Let $z_1=g(w_1)$ and $z_2=g(\zeta_0),$ where
$g(w)=(w-a_j)/\tilde\rho_j.$
Then, by the definition of $D_j$ and Theorem \ref{thm:ct},
$$
D_j(w_1,\zeta_0)=\rho_j D(z_1,z_2)\ge \frac{4\rho_j}{\pi} h_{\D^*}(z_1,z_2).
$$
Since $g$ maps $\Delta_j$ conformally onto $\D^*,$
the principle of the hyperbolic metric leads to the following:
$$
h_{\D^*}(z_1,z_2)=h_{\Delta_j}(w_1,\zeta_0)\ge h_X(w_1,\zeta_0).
$$
When $j=n,$ with $g(w)=\tilde\rho_n/w,$
we have the estimate
$D_n(w_1,\zeta_0)\ge (4\rho_n/\pi)h_X(w_1,\zeta_0)$
in a similar way.
Since $4\rho_j/\pi\ge 2m/\pi=2\min\rho_k/\pi,$ we obtain
$$
d_X(w_1,w_2)\ge N_1\big[ h_X(w_1,\zeta_0)+h_X(\zeta_0,w_2)\big]
\ge N_1h_X(w_1,w_2),
$$
where $N_1=2m/\pi.$
Similarly, we get the first inequality in the other cases with the same constant $N_1.$
Thus the first inequality has been shown.

We next show the inequality $e_X(w_1,w_2)\le N_2h_X(w_1,w_2).$
We consider several cases according to the location of $w_1,w_2.$

\noindent
{\sl Case (i)} $w_1, w_2\in E^*_j:$
We first assume that $j\ne n$ and choose $k\ne j$ so that $\tilde\rho_j=|a_k-a_j|.$
Let $X_1=\C\setminus\{a_j,a_k\}.$
Then $X_1\supset X$ and
$g(w)=(w-a_j)/(a_k-a_j)$ maps $X_1$ conformally onto $\C_{0,1}.$
Set $z_1=g(w_1)$ and $z_2=g(w_2).$
By Theorem \ref{thm:ct}, we obtain
$$
D_j(w_1,w_2)=\rho_j D(z_1,z_2)\le M_0\rho_jh_{\C_{0,1}}(z_1,z_2)
=M_0\rho_jh_{X_1}(w_1,w_2)
\le M_0\rho_jh_X(w_1,w_2),
$$
and thus $e_X(w_1,w_2)\le M_0 \rho_j h_X(w_1,w_2).$
When $j=n,$ we set $g(w)=a_k/w,$ where $a_k$ is chosen so that
$\tilde\rho_n=|a_k|.$
Then, we also have the estimate $D_n(w_1,w_2)\le M_0\rho_n h_X(w_1,w_2).$
In summary, we have $D_j(w_1,w_2)=e_X(w_1,w_2)\le M_0\rho_{n} h_X(w_1,w_2)$
for $j=1,\dots,n,$ because $\rho_j\le e^{-2}\rho_n<\rho_n.$

\noindent
{\sl Case (ii)} $w_1\in E^*_j$ and $w_2\in W:$
Let $\zeta$ be the intersection point of the line segment $[w_1,w_2]$ with
the boundary circle $\partial E_j.$
%Recall that $e_X(w_1,w_2)=D_j(w_1,\zeta_0)+|\zeta_0-w_2|$ by definition.
It is thus enough to show the inequalities
\begin{equation}\label{eq:B}
D_j(w_1,\zeta)\le B_1h_{X}(w_1,w_2)
\aand
|\zeta-w_2|\le B_2h_{X}(w_1,w_2)
\end{equation}
for some constants $B_1$ and $B_2.$

We start with the second one.
Assume that $j\ne n$ for a while.
%Taking a $k$ with $k\ne j,n,$ set $X_1=\C\setminus\{a_j,a_k\}.$
Then the function $g(w)=a_j/w$ maps $X$ conformally into $\C_{0,1}.$
(When $j=1,$ we set $g(w)=1/w.$)
Put $z_l=g(w_l)$ for $l=1,2$ and set $X_1
%=\C\setminus\{g\inv(0), g\inv(1)\}
=g\inv(\C_{0,1})$ which contains $X.$
We consider a shortest hyperbolic geodesic $\alpha$ joining $w_1$ and $w_2$ in $X_1.$
Note that $\ell_{X_1}(\alpha)=h_{X_1}(w_1,w_2)\le h_X(w_1,w_2)$ by
the principle of hyperbolic metric.
In order to complete the proof, we need to analyse the location of
the geodesic $\alpha.$
Since $|w_l|\le\rho_n,~l=1,2,$ the points $z_l$ satisfy
$|z_l|=|a_j/w_l|\ge \sigma,$
where $\sigma=|a_j|/\rho_n\le \tilde\rho_n/\rho_n\le e\inv.$
(When $j=1,$ $\sigma:=1/\rho_n\le e\inv$ because $1\notin X.$)
We now apply Corollary \ref{cor:N} to see that $g(\alpha)$ is contained
in the set $|z|\ge \sigma e^{-K_0},$ where we recall that $K_0\approx 0.85.$
Thus, $\alpha$ lies in the set $|w|\le |a_j|e^{K_0}/\sigma=\rho_ne^{K_0}.$
%We note also that $|a_j|\ge \tilde\rho_1=e\rho_1$ for $1<j<n.$
We note also that $\sigma=|a_j|/\rho_n\ge\tilde\rho_1/\rho_n=e\rho_1/\rho_n.$
Hence, by \eqref{eq:Hempel}, we have the lower estimate
\begin{align*}
\lambda_{X_1}(w)=\lambda(g(w))|g'(w)|
&\ge \frac{1}{2|w|(C_0+|\log|a_j/w||)}\ge \frac{e^{-K_0}}{2\rho_n(C_0-K_0-\log\sigma)} \\
&\ge \frac{e^{-K_0}}{2\rho_n(C_0-K_0-1+\log(\rho_n/\rho_1))}=:\frac1{U_1}
\end{align*}
for $w\in \alpha, 1<j<n.$
Since $e\rho_1\le1,$ this estimate is valid also for $j=1.$

We are now ready to show the second inequality in \eqref{eq:B} for $j\ne n.$
We denote by $\beta$ the line passing through $\zeta$ and orthogonal to the
line segment $[w_1,w_2].$
Take an intersection point $\zeta_1$ of the line $\beta$ with the geodesic $\alpha$
and denote by $\alpha_1$ the part of $\alpha$ joining $\zeta_1$ and $w_2.$
From the above inequality, we derive
$$
h_X(\zeta_1,w_2)
\ge h_{X_1}(\zeta_1,w_2)
=\int_{\alpha_1}\lambda_{X_1}(w)|dw|\ge U_1\inv\int_{\alpha_1}|dw|
\ge U_1\inv|\zeta_1-w_2|.
$$
Now we have
$$
|\zeta-w_2|
\le|\zeta_1-w_2|
\le U_1 h_{X_1}(\zeta_1,w_2)
\le U_1 h_{X_1}(w_1,w_2)
\le U_1 h_{X}(w_1,w_2).
$$

We now turn to the case when $j=n.$
Then we need to modify the above argument a bit.
In this case, we set $g(w)=w$ and $X_1=\C_{0,1}.$
If $\alpha$ is contained in the disk $|w|\le 3\rho_n,$
\eqref{eq:Hempel} yields
\begin{equation}\label{eq:Hn}
\lambda_{X_1}(w)\ge\frac{1}{6\rho_n(C_0+\log 3+\log\rho_n)}=:\frac1{U_2}
\end{equation}
for $w\in\alpha.$
Then the same argument as above yields the inequality $|\zeta-w_2|
\le U_2 h_X(w_1,w_2).$
Otherwise, we define $\zeta_2$ to be the first hitting point of the geodesic
$\alpha$ to the circle $\Gamma=\{w:|w-w_2|=|\zeta-w_2|\}$ from $w_2.$
Let $\alpha_2$ be the part of $\alpha$ joining $\zeta_2$ and $w_2$ as before.
Since the inside of $\Gamma$ is contained in the disk $|w|\le 3\rho_n,$ the inequality
\eqref{eq:Hn} holds for $w\in\alpha_2.$
Thus, we have
$$
|\zeta-w_2|=|\zeta_2-w_2|\le U_2 h_{X_1}(\zeta_2,w_2)
\le U_2 h_{X_1}(w_1,w_2)
\le U_2 h_{X}(w_1,w_2).
$$
In this way, we saw that the second inequality in \eqref{eq:B}
with $B_2=\max\{U_1,U_2\}$ holds at any event.

Next we show the first inequality in \eqref{eq:B}.
By case (i), we have
$$
D_j(w_1,\zeta)\le M_0\rho_n h_X(w_1,\zeta).
$$
On the other hand, by making use of the first part of the theorem and
the second inequality in \eqref{eq:B}, we have
$$
h_X(\zeta,w_2)\le N_1\inv d_X(\zeta,w_2)=N_1\inv|\zeta-w_2|
\le N_1\inv B_2h_X(w_1,w_2).
$$
Combining these inequalities, we get
\begin{align*}
D_j(w_1,\zeta)
&\le M_0\rho_n h_X(w_1,\zeta)
\le M_0\rho_n \big\{h_X(w_1,w_2)+h_X(w_2,\zeta)\big\} \\
&\le M_0\rho_n(1+N_1\inv B_2)h_X(w_1,w_2).
\end{align*}
Thus $D_j(w_1,\zeta)\le B_1h_X(w_1,w_2)$ for $j=1,\dots,n$ with
$B_1=M_0\rho_n(1+N_1\inv B_2).$

We have now shown the inequality $e_X(w_1,w_2)\le N_2' h_X(w_1,w_2)$
in this case, where $N_2'=B_1+B_2.$

\noindent
{\sl Case (iii)} $w_1\in W, w_2\in E_j^*:$
This is essentially same as case (ii).

\noindent
{\sl Case (iv)} $w_1,w_2\in W:$
Similarly, we obtain $|w_1-w_2|\le B_2h_X(w_1,w_2)$
with the same constant $B_2$ as in case (ii).

\noindent
{\sl Case (v)} $w_1\in E_j^*$ and $w_2\in E_k^*$ with $j\ne k:$
Then, by definition,
$$
e_X(w_1,w_2)=D_j(w_1,\zeta_1) +|\zeta_1-\zeta_2|+D_k(\zeta_2, w_2),
$$
where $\zeta_1$ and $\zeta_2$ are the intersection points of the line segment
$[w_1,w_2]$ with $\partial E_j$ and $\partial E_k,$ respectively.
By using the auxiliary lines $\beta_l$ orthogonally intersecting $[w_1,w_2]$
at $\zeta_l~(l=1,2),$ we obtain the inequality
\begin{equation}\label{eq:seg}
|\zeta_1-\zeta_2|\le B_2 h_X(w_1,w_2)
\end{equation}
in the same way as in case (ii).
Let $\alpha$ be a shortest hyperbolic geodesic joining $w_1$ and $w_2$ in $X.$
Let $\zeta_l'$ be an intersection point of $\alpha$ with $\beta_l$ for $l=1,2.$

First assume that $D_j(w_1,\zeta_1)\ge 4\rho_j.$
Since $D_j(\zeta_1,\zeta_1')=|\zeta_1-\zeta_1'|\le 2\rho_j\le D_j(w_1,\zeta_1)/2,$
we obtain $D_j(w_1,\zeta_1)\le 2D_j(w_1,\zeta_1').$
By the first part of the theorem, we now observe that
$$
D_j(w_1,\zeta_1')=d_X(w_1,\zeta_1')\le M_0\rho_{n}h_X(w_1,\zeta_1')
\le M_0\rho_{n}h_X(w_1,w_2).
$$
Thus $D_j(w_1,\zeta_1)\le 2M_0\rho_{n}h_X(w_1,w_2).$
Next assume that $D_j(w_1,\zeta_1)<4\rho_j.$
By \eqref{eq:dist}, we have
$$
|\zeta_1-\zeta_2|\ge \dist(E_j,E_k)\ge(1-e\inv-e^{-2})\rho_j
$$
Thus, with the help of \eqref{eq:seg}, we have
$$
D_j(w_1,\zeta_1)<\frac{4}{1-e\inv-e^{-2}}\,|\zeta_1-\zeta_2|
\le\frac{4B_2}{1-e\inv-e^{-2}}\, h_X(w_1,w_2).
$$
We can deal with $D_k(\zeta_2,w_2)$ in the same way.
Therefore, letting
$$
N_2''=B_2+2\max\left\{2M_0\rho_{n},~ \frac{4B_2}{1-e\inv-e^{-2}}\right\},
$$
we obtain the inequality $e_X(w_1,w_2)\le N_2' h_X(w_1,w_2).$

Summarising all the cases, we now conclude that the right-most inequality
in the assertion of the theorem holds with the choice $N_2=\max\{N_2',N_2''\}.$
\end{pf}

We end the paper with the observation that by the proof we can take the bounds $N_1$ and
$N_2$ in the last theorem under the convention $a_n=\infty$ as follows:
$$
N_1=\frac{2\rho_0}\pi
\aand
N_2=C\, \frac{\rho_n}{\rho_0}\log\frac{\rho_n}{\rho_0},
$$
where $C>0$ is an absolute constant and
$$
\rho_0=\min_{1\le j<n}\rho_j.
$$

\def\cprime{$'$} \def\cprime{$'$} \def\cprime{$'$}
\providecommand{\bysame}{\leavevmode\hbox to3em{\hrulefill}\thinspace}
\providecommand{\MR}{\relax\ifhmode\unskip\space\fi MR }
\providecommand{\MRhref}[2]{%
  \href{http://www.ams.org/mathscinet-getitem?mr=#1}{#2}
}
\providecommand{\href}[2]{#2}

%\bibliography{papers}

%\newpage
\end{document}